\documentclass[10pt]{amsart}
\usepackage{amsmath}
\usepackage{amssymb}
\usepackage{amsthm}
\usepackage{hyperref}
\usepackage{amsmath,amsthm,amssymb}

\DeclareMathOperator{\fracc}{frac}

\DeclareMathOperator{\tr}{tr}
\DeclareMathOperator{\Span}{Span}
\DeclareMathOperator{\dist}{dist}

\newcommand{\Q}{\mathbb{Q}}

\newcommand{\Z}{\mathbb{Z}}
\newcommand{\C}{\mathbb{C}}
\newcommand{\R}{\mathbb{R}}
\newcommand{\HH}{\mathbb{H}}
\newcommand{\m}[4]{\ensuremath{\begin{pmatrix}{#1}&{#2}\\{#3}&{#4}\end{pmatrix}}}

\theoremstyle{plain}
\newtheorem{theorem}{Theorem}[section]
\newtheorem*{theorem*}{Theorem}
\newtheorem{lemma}[theorem]{Lemma}

\newtheorem*{corollary*}{Corollary}

\theoremstyle{remark}

\theoremstyle{definition}
\newtheorem{definition}[theorem]{Definition}
\newtheorem*{definition*}{Definition}
\theoremstyle{remark}
\newtheorem{example}[theorem]{Example}

\title{The Baum-Connes conjecture for countable subgroups of SL(2)}
\author{Dmitry Matsnev}
\thanks{Research supported in part by the Center for Mathematical Analysis, Geometry, and Dynamical Systems (Project FJ11) and by Funda\c{c}\~ao para a Ci\^encia e a Tecnologia through program POCI 2010/FEDER.}
\email{matsnev@math.ist.utl.pt}
\date{\today}
\begin{document}
\address{Departamento de Matem\'atica, Instituto Superior T\'ecnico, Av. Rovisco Pais, 1049-001 Lisboa, Portugal}
\maketitle
\begin{abstract}
We present an alternative approach to the result of Guentner, Higson, and Weinberger concerning the Baum-Connes conjecture for finitely generated subgroups of $SL(2,\C)$. Using finite-dimensional methods, we show that the Baum-Connes assembly map for such groups is an isomorphism.
\end{abstract}
\section{Introduction}\label{chapterIntroduction}

The Baum-Connes conjecture, introduced in the early 80's by Paul Baum and Alain Connes, connects the $K$-theory of the reduced crossed product of a $C^*$-algebra by a group acting on such algebra and the $K$-homology of the corresponding classifying space of proper actions of that group (for a formal account see~\cite{BaumConnesHigson}).

Let $\Gamma$ be a discrete group acting on a $C^*$-algebra $A$ by automorphisms. The Baum-Connes conjecture proposes that the ``assembly''  morphism 
$$\mu : KK^{\Gamma}(\underline{E}\Gamma,A)\to K(A\rtimes_r\Gamma)$$
from the $K$-homology of the classifying space $\underline{E}\Gamma$ of proper actions of $\Gamma$ to the $K$-theory of the reduced crossed product of $A$ by $\Gamma$ is an isomorphism.

While the conjecture is formulated in terms of pairs $(\Gamma,A)$, it is possible to state it purely in terms of group $\Gamma$: one can ask whether the original conjecture holds for the group $\Gamma$ and any $C^*$-algebra $A$ on which $\Gamma$ acts. This version of the conjecture is called \emph{the Baum-Connes conjecture with coefficients}\footnote{Some authors refer to it as to the Baum-Connes \emph{property} with coefficients, in the view of the counterexamples (modulo a statement due to Gromov) by Higson, Lafforgue, and Skandalis in~\cite{HigsonLafforgueSkandalis}.}, and this will be our main concern in this work.

To indicate some connections of this conjecture with other areas, we mention that the injectivity of the assembly map implies the Novikov's higher signature conjecture, while the surjectivity of the assembly map has to do with the Idempotents conjecture of Kadison and Kaplansky. A more comprehensive account on various versions of the Baum-Connes conjecture and the ambient areas consult~\cite{MislinValette}.

The starting point of this work is the following 

\begin{theorem}[Guentner, Higson, and Weinberger,~\cite{GuentnerHigsonWeinberger}]\label{theoremGL2Std}
For any field $K$, the Baum-Connes conjecture with coefficients holds for any countable subgroup of $GL(2,K)$.
\end{theorem}

The original proof of this theorem was based on the existence of a metrically proper and  isometric action of such group on a Hilbert space and then appealing to the result of Higson and Kasparov in~\cite{HigsonKasparov}. In this construction the group action was not constructed explicitly, and the Hilbert space on which the groups act can easily be infinite-dimensional. We are interested in giving a more elementary proof of Theorem~\ref{theoremGL2Std}, without appealing to Hilbert space techniques and with a more direct description of the group action.

As a result, we prove the following

\begin{theorem}\label{theoremMain2}
Let $K$ be any field of characteristic 0, and $\Gamma$ be a finitely generated subgroup of $SL(2,K)$. Then $\Gamma$ satisfies the Baum-Connes conjecture with coefficients.
\end{theorem}

Of course, one can also extend this result to any countable subgroup $\Gamma$ of $SL(2,K)$, since $\Gamma$ is a directed union of its finitely generated subgroups, and the Baum-Connes conjecture holds for directed unions of groups satisfying the conjecture (see~\cite{MislinValette}).

\section{Technical tools}\label{sectionTechnicalTools}

In this section we collect some technical tools to be used later in the discussion.

\subsection{Discrete valuations}\label{sectionPreliminariesDiscreteValuations}

\begin{definition}\label{definitionDiscreteValuation}
Let $R$ be an integral domain (a commutative ring without zero divisors in which $0\ne1$). A map $\nu:R\to\Z\cup\{+\infty\}$ is called a \emph{discrete valuation} if it satisfies the following properties for any $a, b\in R$:
\begin{itemize}
\item $\nu(a)=+\infty$ if and only if $a=0$
\item $\nu(ab)=\nu(a)+\nu(b)$
\item $\nu(a+b)\geq\min(\nu(a),\nu(b))$
\end{itemize}
\end{definition}

Given a discrete valuation $\nu$ of an integral domain $R$, it can be extended to the field of fractions $\fracc(R)$ by
\begin{equation*}
\nu\left(\frac ab\right)=\nu(a)-\nu(b),\qquad\qquad a, b\in R.
\end{equation*}

Starting from a field $K$, one can regard $K$ as a ring and define the corresponding notion of a discrete valuation, following Definition~\ref{definitionDiscreteValuation}. Notice that the extension to the field of fractions is consistent with such treatment.

To any discrete valuation $\nu$ one associates its \emph{ring of integers} $\mathcal O_{\nu}$. It consists of all elements of the field $K$ with non-negative valuation. Any element $\pi$ of $\mathcal O_{\nu}$ with $\nu(\pi)=1$ is called a \emph{uniformizer} of $\nu$.

\begin{example}
For any prime number $p$ the $p$-adic valuation $\nu_p$ on $\Q$ is defined by
\[\nu_p\left(p^n\frac ab\right)=n, a, b, n\in\mathbb Z \mbox{ and } (a,p)=(b,p)=1.\]
The ring of integers of $\nu_p$ consists of all rational numbers without any occurrence of $p$ in the denominator, $p$ itself serves as a uniformizer, and the residue field $\mathcal O_{\nu}/\pi\mathcal O_{\nu}$ is isomorphic to $\Z/p\Z$.
\end{example}

\begin{example}
Another construction which we shall employ is an extension of the $p$-adic valuation $\nu_p$ to an algebraic field extension of $\Q$. Suppose that such an extension is $\Q(\gamma)$, where $\gamma$ has degree $m$. Then the extension of $\nu_p$ is defined as
\[\widetilde{\nu_p}(q_{m-1}\gamma^{m-1}+q_{m-2}\gamma^{m-2}+\dots+q_0)=\min\{\nu_p(q_{m-1}),\nu_k(q_{m-2}),\dots,\nu_k(q_0)\}.\]
(here $q_{m-1},q_{m-2},\dots,q_0\in\Q$.)
\end{example}

\subsection{Simplicial tree for $SL(2,K)$}\label{subsectionSimplicialTree}

Suppose $K$ is a field with a discrete valuation $\nu$, its ring of integers $\mathcal O$, and a uniformizer $\pi$.

A simplicial tree $T_{\nu}$ of equivalence classes of lattices in $K^2$ is constructed as follows. The vertices of the tree are the homothety classes of $\mathcal O_{\nu}$-lattices in $K^2$. Two such vertices are connected by an edge if there exist some lattice representatives $L$ and $M$ of these vertices such that
\[\pi L\subset M\subset L.\]

Notice that the valence of each vertex of $T_{\nu}$ is equal to the cardinality of the residue field.

The action of $SL(2,K)$ on $T_{\nu}$ is defined via a natural action on the corresponding lattices. If we equip the tree with a standard simplicial metric then the action is an isometry. For more in-depth discussion of this construction, consult~\cite{Serre}.

\subsection{Groups of integral characteristic and Alperin-Shalen reduction}

In their study of the cohomological dimension of linear groups in~\cite{AlperinShalen}, Alperin and Shalen introduced a reduction technique from a given linear group to a family of its subgroups of a special type, which we shall summarize here.

\begin{definition}
A subgroup $\Gamma$ of $SL(2,\C)$ is said to have an \emph{integral characteristic} if the coefficients of the characteristic polynomial of every element of $\Gamma$ are algebraic integers. 
\end{definition}

We shall rely on the following technical fact proven in~\cite{AlperinShalen}.

\begin{theorem}\label{theoremAlperinShalen}
Let $\Gamma$ be a finitely generated subgroup of $SL(2, \C)$ and let $A$ be a (finitely generated) ring of the matrix entries of the elements of $\Gamma$ with its field of fractions being $K$. Then there exist finitely many discrete valuations $\nu_1, \nu_2, \dots, \nu_m$ on $K$ and a finite sequence (``hierarchy'')
\[\mathcal H_0, \mathcal H_1, \dots, \mathcal H_m\]
of families of subgroups of $\Gamma$ such that $\mathcal H_m$ consists of $\Gamma$ only, $\mathcal H_0$ consists of some subgroups of $\Gamma$ of integral characteristic, and  for each $i$ with $1\leq i\leq m$ any group $G$ from $\mathcal H_i$ acts on a simplicial tree of $\mathcal O_{\nu_i}$-equivalence classes with the isometry of the action belonging to $\mathcal H_{i-1}$.
\end{theorem}

\subsection{Metrically proper actions}

\begin{definition}
Let $\Gamma$ be a discrete group acting on a metric space $X$ by isometries. The action is called \emph{metrically proper} if for every bounded subset $B$ of $X$ the set
\[\{g\in\Gamma | g.B\cap B\ne\emptyset \}\]
is finite.
\end{definition}

We remark here that the condition
\[\forall x\in X \; \forall C \; \{g\in\Gamma | \dist(x,g.x)<C \} \mbox{ is finite}\]
implies that the action is metrically proper (to see this, one can take $C$ to be three times the diameter of $B$ and use the fact that the action is an isometry.)

\section{Proof of the theorem}

In what follows, we shall extensively use the fact that the Baum-Connes conjecture with coefficients passes to subgroups (see~\cite{ChabertEchterhoff}), without explicitly mentioning it.

The argument will be structured as follows. In Section~\ref{sectionSubgroupsAlgebraic} we prepare the ground by constructing a metrically proper action for finitely generated subgroups over an algebraic extension of $\Q$.

In Section~\ref{sectionReductionToIntegral} we prove Theorem~\ref{theoremReductionToIntChar}, which shows that in order to prove the Baum-Connes conjecture for a given finitely generated subgroup of $SL(2,\C)$, it is sufficient to prove it for the subgroups of integral characteristic of the original group. 

In Section~\ref{sectionZariskiDense} we treat the case of subgroups of integral characteristic which are Zariski-dense in $SL(2,\C)$. Namely, we prove in  Theorem~\ref{theoremZariskiDense} that for all such subgroups coming from a finitely generated group in $SL(2,\C)$ we started with, the Baum-Connes conjecture holds by means of using the results from Section~\ref{sectionSubgroupsAlgebraic}.

Finally, in Section~\ref{sectionZariskiNondense} we prove in  Theorem~\ref{theoremZariskiNondense} that all countable subgroups of $SL(2,\C)$ of integral characteristic which are not Zariski-dense in it satisfy the Baum-Conjecture.

\subsection{Subgroups over algebraic fields}\label{sectionSubgroupsAlgebraic}

We start with a technical construction of a certain proper action which will be used later.

\begin{lemma}\label{lemmaProperActionAlgebraic}
Let $\Gamma$ be a finitely generated subgroup of $SL(2,K)$, where $K$ is a finite extension of $\Q$. Then $\Gamma$ acts metrically properly on a finite-dimensional space (a finite product of simplicial trees and a hyperbolic plane). 
\end{lemma}
\begin{proof}
We can treat $\Gamma$ as defined over a finitely generated ring which we can take, enlarging it if necessary, to be $A=\Z[\frac1s,\gamma]$, where $s$ is a natural number (the l.u.b. of all the denominators of the rationals participating in the entries of the generating set) and $\gamma$ is an algebraic integer. We assume that prime factorization of $s$ is $p_1\cdot p_2\cdots p_n$ with nonrepeating terms.

For each $p_k$, $k=1,2,\dots,n$ consider the $p_k$-adic valuation $\nu_k$ on $\Q\cap A$ and extend it to the whole $A$ and its fraction field by the usual rule.
Denote by $T_k$ the simplicial tree corresponding to $\widetilde{\nu_k}$ and denote by $\alpha_k$ the induced action of $\Gamma$ on this tree.

Define an action $\alpha_H$ of $\Gamma$ on the 2-dimensional real hyperbolic space $\HH_2$ via the natural isometric action of $SL(2,\R)$ on $\HH_2$.

We claim that the diagonal action 
\[\alpha=\alpha_1\times\alpha_2\times\dots\times\alpha_n\times\alpha_H\]
on the product of trees and a hyperbolic space $\HH_2$
\[T=T_1\times T_2\times\dots\times T_n\times\HH_2\]
is metrically proper. To clarify this, we shall show that for any bounded set $B$
\[\#\{g\in\Gamma|g.B\cap B\neq\emptyset\}<\infty.\]

Since it is enough to prove the statement for sufficiently large sets $B$ only, we shall enlarge $B$ in the following way: let $B_k$ be the projection of $B$ to the $k$-th tree, and $B_{\mathbb H_2}$ be the projection of $B$ to the hyperbolic space. Clearly $B\subseteq B_1\times\dots\times B_n\times B_{\mathbb H_2}$,
and we shall be working with the latter set instead of the original $B$.

Now the set $\{g\in\Gamma|g.B\cap B)\neq\emptyset\}$ is actually
\[\bigcap_{k=1}^n\{g\in\Gamma|g.B_k\cap B_k\neq\emptyset\}\cap\{g\in\Gamma|g.B_{\mathbb H_2}\cap B_{\mathbb H_2}\neq\emptyset\},\]
therefore it is enough to show that for any number $C$ and any points $v_k\in T_k$ ($k=1,\dots, n$) and $x\in\mathbb{H}_2$,
\begin{multline}\label{propernessWithC}
\#\{g\in\Gamma|\dist_{T_k}(g.v_k,v_k)<C \;\mbox{for all }k=1,\dots,n,\\\;\mbox{and}\;\dist_{\mathbb{H}_2}(g.x,x)<C\}<\infty.
\end{multline}

Further, by triangle inequality it is sufficient to check condition (\ref{propernessWithC}) for the ``root'' vertices $v_0$ of the trees and the ``center'' point $x_0$ of the hyperbolic space.

Now let $g=(g_{ij})\in\Gamma$ be an element from the set in question. Here each matrix entry $g_{ij}$ belongs to $\Z[\frac1{p_1},\dots, \frac1{p_n},\gamma]$, say
\[g_{ij}=\sum_{0\leq l\leq m-1}\frac{a_{ijl}}{\prod_{k=1}^np_k^{n_{ijkl}}}\gamma^l,\qquad a_{ijl}, n_{ijkl}\in\Z,\; (a_{ij},p_k)=1, \; k=1,\dots, n.\]

The distance from $g.v_0$ to $v_0$ in the tree $T_{p_k}$ is bounded only if each matrix entry has limited from above powers $n_{ijkl}$ in the denominator, which means the denominator itself should be bounded.

Since $\dist_{\mathbb{H}_2}(g.x_0,x_0)<C$ means $\cosh(\dist_{\mathbb{H}_2}(g.x_0,x_0))<\cosh C$, and 
$$\cosh(\dist_{\mathbb{H}_2}(g.x_0,x_0))=\sum g_{ij}^2,$$
each matrix entry should be bounded. Together with the previous observation, this leads to a finite number of choices for $a_{ijl}$ and $n_{ijkl}$, and thus there are finitely many such elements $g$ in the intersection.
\end{proof}

\subsection{Reduction to groups of integral characteristic}\label{sectionReductionToIntegral}

In this section we applythe Alperin-Shalen hierarchy construction which will allow us to reduce the proof of the Baum-Connes conjecture for any finitely generated subgroup of $SL(2,\C)$ to subgroups of integral characteristic.

The main motivation for the study of the isotropy of group actions on trees is the following result:

\begin{theorem}[Oyono-Oyono,~\cite{OyonoOyono}]\label{theoremOyonoOyono}
Let $\Gamma$ be a discrete countable group acting on a tree. Then the Baum-Connes conjecture holds for $\Gamma$ if and only of it holds for all the isotropy subgroups of the action on vertices of the tree.
\end{theorem}

\begin{theorem}\label{theoremReductionToIntChar}
Let $\Gamma$ be a finitely generated subgroup of $SL(2,\C)$. Then the Baum-Connes conjecture with coefficients holds for $\Gamma$ if and only if it holds for all subgroups of $\Gamma$ of integral characteristic.
\end{theorem}
\begin{proof}
Given $\Gamma$, apply Theorem~\ref{theoremAlperinShalen} and consider a hierarchy of families of subgroups of $\Gamma$ together with the actions on trees which that theorem furnishes. Repeated applications of Theorem~\ref{theoremOyonoOyono} allow one to reduce the Baum-Connes conjecture for the top level of the hierarchy (for $\Gamma$ itself, that is) to the one for the bottom of the hierarchy, which contains subgroups of $\Gamma$ of integral characteristic only.
\end{proof}

\subsection{Zariski-dense subgroups}\label{sectionZariskiDense}

Now we need to prove the Baum-Connes conjecture for subgroups of integral characteristic. In this section we concentrate on integral characteristic subgroups $\Gamma$ of $SL(2,\C)$ whose Zariski closure is the entire $SL(2,\C)$.
 
The following result essentially goes back to Zimmer (cf.~\cite[Lemma~6.1.7]{ZimmerErgodicTheory}).

\begin{lemma}\label{lemmaZimmerEmbedding}
Let $\Gamma$ be a Zariski-dense subgroup of $SL(2,\C)$ of integral characteristic. Then there exists a faithful representation
\[\alpha: SL(2,\C)\to GL(4,\C)\]
such that the matrix entries of every element of $\alpha(\Gamma)$ is an algebraic integer.
\end{lemma}
\begin{proof}
We shall write $G$ for the Zariski-closure of $\Gamma$ in $SL(2,\C)$, that is, $G=SL(2,\C)$.  Let $f_g$ be a map $G\to\C$ defined by
\[f_g:h\mapsto\tr(gh), \qquad h\in G.\]
Notice that $f_{g_1}(h)+f_{g_2}(h)=\tr((g_1+g_2)h)$ and $\lambda f_g(h)=\tr((\lambda g)h)$ for any $g_1, g_2, h\in G$ and $\lambda\in\C$. This allows us to consider $f_g(h)$ as a short-hand notation for $\tr(gh)$ for any $g\in\C G$ and $h\in G$. Let
\[V=\Span_{\C}\langle f_g \rangle_{g\in\C G}.\]
Since the conditions defining $f_g$ are linear with respect to the entries of $g$, the linear space $V$ has finite dimension (more precisely, its dimension is $4$.)

Consider the following action of $G$ on $V$:
\begin{equation}\label{equationSL2ActionZDense}
g.f_h=f_{gh}, \qquad g\in G, h\in\C G.
\end{equation}
Since this action is linear, we have a representation of $G$.

Let
\[W=\Span_{\C}\langle f_g \rangle_{g\in\C\Gamma}.\]
This subspace is $\Gamma$-invariant and, since $\Gamma$ is Zariski-dense in $G$, is also $G$-invariant. Thus $W=V$.

Let $g_1, g_2, g_3, g_4\in\Gamma$ be such that $\left\lbrace f_{g_1}, f_{g_2}, f_{g_3}, f_{g_4}\right\rbrace $ is a basis of $V$ (we can arrange this because $V$ is generated by $f_g$ for $g\in\Gamma$). With respect to this basis the action~\eqref{equationSL2ActionZDense} is given by a matrix $(\alpha_{ij}^g)$, such that
\begin{equation}\label{equationSL2ActionMatrix}
g.f_{g_i}(h)=\sum_{j=1}^4\alpha_{ij}^gf_{g_j}(h), \qquad g, h\in G, i=1,\dots,4.
\end{equation}
Thus we obtain a representation $\alpha: G\to GL(4,\C)$.

We confirm that $\alpha$ is faithful by showing that from the identity
$g.f_1(h)=f_1(h)$ for all $h\in G$ follows that $g=1$. To do this, take elementary matrices for $h$ and write this condition entries-wise.

If $g\in\Gamma$, \eqref{equationSL2ActionMatrix} means that in particular
\[\tr(gg_ig_k)=g.f_{g_i}(g_k)=\sum_{j=1}^4\alpha_{ij}^gf_{g_j}(g_k)=\sum_{j=1}^4\alpha_{ij}^g\tr(g_jg_k), \qquad i, k=1,\dots,4.\]
Then $\{\alpha_{ij}^g\}$ are the solutions of a system of linear equations with algebraic coefficients, and therefore the matrix entries $(\alpha_{ij}^g)$ of the representation $\alpha$ are algebraic.
\end{proof}

\begin{theorem}\label{theoremZariskiDense}
Let $G$ be a finitely generated subgroup of $SL(2,\C)$. Then the Baum-Connes conjecture with coefficients holds for any subgroup $\Gamma$ of $G$, provided that $\Gamma$ has integral characteristic and is Zariski-dense in $SL(2,\C)$. 
\end{theorem}
\begin{proof}
We may assume that given $G$ contains at least one subgroup $\Gamma$ which is both of integral characteristic and is Zariski-dense. Select one such $\Gamma$ and apply  Lemma~\ref{lemmaZimmerEmbedding}, yielding an embedding $\alpha: SL(2,\C)\to GL(4,\C)$. Notice that if $\Gamma'$ is any Zariski-dense subgroup of $G$ of integral characteristic then its image $\alpha(\Gamma')$ is conjugate to a subgroup whose matrix entries are algebraic. Moreover, since $\alpha(G)$ is finitely generated, all matrix entries of its elements belong to a finitely generated subring of $\C$, and this is true for its subgroup $\alpha(\Gamma')$ as well, which means that the entries of the elements of $\alpha(\Gamma')$ belong to a certain finitely generated subring of algebraic numbers which depends on the original $G$ only, rather than on $\Gamma'$.

By the Primitive Element Theorem $\alpha(\Gamma')\subseteq GL(4,K)$ for some field $K$  with $[K:\Q]<\infty$. Take $H=\alpha(SL(2,\C))$. We see that both $\alpha(\Gamma')$ and $H\cap GL(4,K)$ are Zariski-dense in $H$, thus they are both defined over $K$ by~\cite[Proposition 3.1.8]{ZimmerErgodicTheory}. This means $\alpha(\Gamma')$ is locally isomorphic to a subgroup of $SL(2,K)$ (see~\cite[Theorem~7]{Zimmer}). We can represent $\alpha(\Gamma')$ as a directed union of its finitely generated subgroups and for each such subgroup apply Lemma~\ref{lemmaProperActionAlgebraic}, thus obtaining a metrically proper action of on a finite-dimensional space. Via a finite-dimensional version of the result of Higson and Kasparov in~\cite{HigsonKasparov}, the Baum-Connes conjecture for such subgroup follows, hence it follows for $\Gamma'$ as well. 
\end{proof}

\subsection{Zariski-non-dense subgroups}\label{sectionZariskiNondense}

Finally, we discuss the case of Zariski-non-dense subgroups of integral characteristic. The rest of this section is devoted to the proof of the following
\begin{theorem}\label{theoremZariskiNondense}
Let $\Gamma$ be a countable non-Zariski-dense subgroup of $SL(2,\C)$ of integral characteristic. Then the Baum-Connes conjecture with coefficients holds for $\Gamma$.
\end{theorem}

Let us start with some preliminary remarks on algebraic Lie groups.
Suppose $G$ is a Zariski--closed proper subgroup of $SL(2,\C)$. We write $G_0$ for the Zariski-connected component of the unit of $G$. It is known that $G_0$ is a normal subgroup of $G$ of finite index \cite[I.1.2]{Borel}.

Since for algebraic groups the notions of connected and irreducible components coincide~\cite[AG.17.2]{Borel}, $G_0$ is abelian if and only if its Lie algebra is commutative \cite[IV.4.3]{Hochschild}. 
Since $G$ is a proper subgroup of $SL(2,\C)$, its dimension is strictly less than $3$. In the subsections below we shall address each dimension case separately.

\subsection*{Dimension $0$}

In this case $\dim G_0=0$ as well, and, since $G_0$ is connected, we conclude that it is trivial. The group $G$ itself, being a finite extension of $G_0$, is finite, whence the Baum-Connes conjecture for $G$ holds trivially.

\subsection*{Dimension $1$}

The Lie algebra of $G$ (and $G_0$ as well) has to be $1$-dimensional. In particular, it has to be commutative, thus $G_0$ is abelian. There are only two (up to conjugacy) connected abelian $1$-dimensional groups, namely $\left\{\left. \m{a}{0}{0}{a^{-1}}\right| a\in\C^{\times}\right\}$ and $\left\{\left. \m{1}{b}{0}{1}\right| b\in\C\right\}$. We shall treat them separately.

Suppose $G_0=\left\{\m{a}{0}{0}{a^{-1}}\right\}$ (up to conjugacy). Then, since $G_0$ is normal in $G$, any conjugate of $\m{a}{0}{0}{a^{-1}}$ by any element in $G$, say $\m{g_{11}}{g_{12}}{g_{21}}{g_{22}}$, has to have the same diagonal form:
\begin{multline*}
\m{g_{11}}{g_{12}}{g_{21}}{g_{22}}\m{a}{0}{0}{a^{-1}}\m{g_{11}}{g_{12}}{g_{21}}{g_{22}}^{-1}=\\
\m{*
}{(a^{-1}-a)g_{11}g_{12}}{(a-a^{-1})g_{21}g_{22}}{*
}=
\m{b}{0}{0}{b^{-1}}.
\end{multline*}
This means that $g_{11}g_{12}=0$ and $g_{21}g_{22}=0$. To satisfy the first condition, we need to take either $g_{11}$ to be zero or $g_{12}$ to be zero. Thus $G$ may contain only matrices with zeros on the diagonal, or off the diagonal:
\[G\subseteq\left\{\m{a}{0}{0}{a^{-1}}, \m{0}{a}{-a^{-1}}{0}\right\}=H.\]
This group $H$ is amenable, and, modifying Theorem~\ref{theoremOyonoOyono}, it is possible to show (see~\cite[Theorems 5.18 and 5.23]{MislinValette}) that any countable subgroup of $H$ satisfies the Baum-Connes conjecture with coefficients\footnote{Theorem 5.18 in~\cite{MislinValette} shows that $H$ satisfies the Baum-Connes conjecture with trivial coefficients, while Theorem 5.23 proves that any countable subgroup of such group satisfies the conjecture with arbitrary coefficients.}. 

Now suppose $G_0=\left\{\m{1}{b}{0}{1}\right\}$ (again, up to conjugacy). Since $G_0$ is normal in $G$, any conjugate of $\m{1}{b}{0}{1}$ by an arbitrary element $\m{g_{11}}{g_{12}}{g_{21}}{g_{22}}$ element of $G$ should have the same form:
\[\m{g_{11}}{g_{12}}{g_{21}}{g_{22}}\m{1}{b}{0}{1}\m{g_{11}}{g_{12}}{g_{21}}{g_{22}}^{-1}=\m{*}{*}{-bg_{21}^2}{*}=\m{*}{*}{0}{*}.\]
Thus we have $g_{21}=0$, which means that $G$ contains only matrices of the form $\m{a}{b}{0}{a^{-1}}$, and we shall discuss this group in the following subsection.

\subsection*{Dimension $2$}

Let $H$ denote the subgroup of $SL(2,\mathbb C)$ consisting of all the matrices of the form $\m{a}{b}{0}{a^{-1}}$, where $a\in\C^{\times}, b\in\C$. Note that its Lie algebra consists of the matrices of the form $\m{a}{b}{0}{-a}$, $a, b\in\C$.
\begin{lemma}\label{MaximalityOfH}
Any subgroup $K$ of $SL(2,\mathbb C)$, which includes $H$ and some element not in $H$, coincides with the whole group  $SL(2,\mathbb C)$.
\end{lemma}
\begin{proof}
Suppose $K$ contains some element \m{g_{11}}{g_{12}}{g_{21}}{g_{22}} with $g_{21}\ne0$. Then we can multiply this element by an element \m{g_{21}^{-1}}{-g_{22}}{0}{g_{21}} of $H$ on the right to get \m{g_{11}g_{21}^{-1}}{-1}{1}{0}.

Now for any complex numbers $a$ and $b$ with $a\ne0$ we can multiply \m{g_{11}g_{21}^{-1}}{-1}{1}{0} on the left by \m{a}{b-ag_{11}g_{21}^{-1}}{0}{a^{-1}} to get \m{b}{-a}{a^{-1}}{0}. Since $K$ is a group, it ought to contain all inverses as well, in particular \m{0}{a}{-a^{-1}}{b}. This means that $K$ contains all matrices of determinant $1$ with $0$ in the upper left corner.

Finally, let us take an arbitrary element of $SL(2,\mathbb C)$, say \m{s_{11}}{s_{12}}{s_{21}}{s_{22}}. Since we already know that all matrices with  $s_{21}=0$ belong to $H$, and therefore to $K$, the essential part of the argument is to show that any such matrix with $s_{21}\ne0$ belongs to $K$. The identity
\[\m{s_{11}}{s_{12}}{s_{21}}{s_{22}}=\m{1}{s_{11}s_{21}^{-1}}{0}{1}\m{0}{-s_{21}^{-1}}{s_{21}}{s_{22}}\]
completes the proof, since all matrices on the right-hand-side belong to $K$.
\end{proof}

Now we provide some technical results about Lie subalgebras of $\mathfrak{sl}(2,\C)$.
\begin{lemma}
For any 2-dimensional noncommutative Lie algebra there exists a basis $\{x_1, x_2\}$ with multiplication rule $[x_1x_2]=x_1$.
\end{lemma}

\begin{lemma}\label{OnlyOne}
The Lie algebra $\mathfrak{sl}(2,\C)$ contains only one (up to conjugation) $2$-dimensional Lie subalgebra, namely $\left\{\left.\m{a}{b}{0}{-a}\right| a, b\in\C\right\}$.
\end{lemma}
\begin{proof}
Suppose we have a $2$-dimensional noncommutative Lie subalgebra $\mathfrak h$ of $\mathfrak{sl}(2,\C)$. Let $\{x_1, x_2\}$ denote the basis of $\mathfrak h$ constructed in the lemma above. A priori there could be two possibilities: both eigenvalues of the matrix $x_2$ coincide (and therefore are zeros) or they are distinct. In the first case $x_2$ is conjugate to its Jordan form, namely $\m{0}{1}{0}{0}$, and if $x_1$ after same conjugation has the form $\m{a}{b}{c}{d}$, then the multiplication condition $[x_1x_2]=x_1$ is
\[\left[\m{a}{b}{c}{d},\m{0}{1}{0}{0}\right]=\m{-c}{a-b}{0}{c}=\m{a}{b}{c}{d},\]
from where we conclude that $a=b=c=d=0$, which means $x_1=0$ and therefore can not serve as basis element, so that the case where both eigenvalues of the matrix $x_2$ coincide can not happen. Now suppose that the eigenvalues of $x_2$ are distinct, say $\lambda$ and $-\lambda$. Conjugating $x_1$ and $x_2$, we write the multiplication condition as
\[\left[\m{a}{b}{c}{d},\m{\lambda}{0}{0}{-\lambda}\right]=\m{0}{-2b\lambda}{2c\lambda}{0}=\m{a}{b}{c}{d},\]
so that we have $a=d=0$ and $2c\lambda=c$, $-2b\lambda=b$. We are looking for solutions with at least one of the coefficients $b$ and $c$ being non-zero, therefore we end up with two possibilities:
\begin{enumerate}
\item $b=0\ne c$, $\lambda=\frac12$
\item $c=0\ne b$, $\lambda=-\frac12$
\end{enumerate}
Thus any non-commutative $2$-dimensional Lie subalgebra of $\mathfrak{sl}(2,\mathbb C)$ is conjugate-equivalent to $\mathfrak h_1=\mathbb C\m{0}{0}{1}{0}\oplus\mathbb C\m{\frac12}{0}{0}{-\frac12}$ or $\mathfrak h_2=\mathbb C\m{0}{1}{0}{0}\oplus\mathbb C\m{-\frac12}{0}{0}{\frac12}$. Finally, $\mathfrak h_1$ and $\mathfrak h_2$ are conjugate to each other via the matrix \m{0}{1}{1}{0}. By scaling the second matrix in $\mathfrak h_2$ , we obtain the representation $\left\{\left. \m{a}{b}{0}{-a} \right| a, b\in\mathbb C\right\}$.

Now we show that $\mathfrak{sl}(2,\C)$ does not contain any commutative $2$-dimensional subalgebras. Suppose one such exists and has a basis $\{x, y\}$. Conjugating by some matrix, we can put $y$ into Jordan form, and let $x$ be represented by \m{x_{11}}{x_{12}}{x_{21}}{-x_{11}} under the same conjugation. We have two possibilities: the eigenvalues of the matrix, representing $y$, coincide (and therefore are zeros) or they are distinct, and by scaling the matrix we assume that they are $1$ and $-1$. In the first case the commutativity condition can be written as
\[\m{x_{11}}{x_{12}}{x_{21}}{-x_{11}}\m{0}{1}{0}{0}=\m{0}{1}{0}{0}\m{x_{11}}{x_{12}}{x_{21}}{-x_{11}},\]
which leads to $x_{11}=x_{21}=0$, so that $x$ is a scalar multiple of $y$, and this can not happen. In the second case we have
\[\m{x_{11}}{x_{12}}{x_{21}}{-x_{11}}\m{1}{0}{0}{-1}=\m{1}{0}{0}{-1}\m{x_{11}}{x_{12}}{x_{21}}{-x_{11}},\]
this means $x_{12}=x_{21}=0$, and again we have a contradiction with linear independence of $x$ and $y$.
\end{proof}

Getting back to the group $G_0$, we see that Lemma~\ref{OnlyOne} describes the Lie algebra of $G_0$, up to conjugacy. Therefore $G_0$ and $H$ are conjugate to each other. 
Lemma~\ref{MaximalityOfH} confirms that there are no proper subgroups of $SL(2,\C)$ larger than $H$, therefore we conclude that $G=G_0$.

Finally, $G$ is a semidirect product
\[\left\{\left. \begin{pmatrix}1&b\\0&1\end{pmatrix} \right| b\in\C\right\}\rtimes\left\{\left. \begin{pmatrix}a&0\\0&a^{-1}\end{pmatrix} \right| a\in\C^{\times}\right\}\]
of two abelian groups, hence it is amenable, and we conclude that any finitely generated subgroup of $G$ (and hence any countable one) satisfies the Baum-Connes conjecture with coefficients by applying~\cite[Theorem 5.23]{MislinValette}.

Now Theorem~\ref{theoremMain2} has been proven completely.

\section*{Acknowledgements}This paper is based on one part of the author's thesis completed under the supervision of Nigel Higson at Penn State University. The author is very grateful to Professor Higson for his invaluable guidance, comments, and suggestions. Also the author is thankful to an anonymous referee for her or his comments and suggestions.
\bibliographystyle{amsalpha}
\addcontentsline{toc}{chapter}{Bibliography}
\bibliography{bibliography-data}
\end{document}